\documentstyle[12pt]{article}
\def\hybrid{\topmargin 0pt      \oddsidemargin 0pt
        \headheight 0pt \headsep 0pt
        \textwidth 6.35in       % A4 paper
        \textheight 9.5in       % A4 paper
        \marginparwidth 0.0in
        \parskip 5pt plus 1pt   \jot = 1.5ex}
\catcode`\@=11
\def\marginnote#1{}

\newcount\hour
\newcount\minute
\newtoks\amorpm
\hour=\time\divide\hour by60
\minute=\time{\multiply\hour by60 \global\advance\minute by-\hour}
\edef\standardtime{{\ifnum\hour<12 \global\amorpm={am}%
        \else\global\amorpm={pm}\advance\hour by-12 \fi
        \ifnum\hour=0 \hour=12 \fi
        \number\hour:\ifnum\minute<10 0\fi\number\minute\the\amorpm}}
\edef\militarytime{\number\hour:\ifnum\minute<10 0\fi\number\minute}

\def\draftlabel#1{{\@bsphack\if@filesw {\let\thepage\relax
   \xdef\@gtempa{\write\@auxout{\string
      \newlabel{#1}{{\@currentlabel}{\thepage}}}}}\@gtempa
   \if@nobreak \ifvmode\nobreak\fi\fi\fi\@esphack}
        \gdef\@eqnlabel{#1}}
\def\@eqnlabel{}
\def\@vacuum{}
\def\draftmarginnote#1{\marginpar{\raggedright\scriptsize\tt#1}}

\def\draftlabel#1{{\@bsphack\if@filesw {\let\thepage\relax
   \xdef\@gtempa{\write\@auxout{\string
      \newlabel{#1}{{\@currentlabel}{\thepage}}}}}\@gtempa
   \if@nobreak \ifvmode\nobreak\fi\fi\fi\@esphack}
        \gdef\@eqnlabel{#1}}
\def\@eqnlabel{}
\def\@vacuum{}
\def\draftmarginnote#1{\marginpar{\raggedright\scriptsize\tt#1}}

\def\draft{\oddsidemargin -.5truein
        \def\@oddfoot{\sl preliminary draft \hfil
        \rm\thepage\hfil\sl\today\quad\militarytime}
        \let\@evenfoot\@oddfoot \overfullrule 3pt
        \let\label=\draftlabel
        \let\marginnote=\draftmarginnote
   \def\@eqnnum{(\theequation)\rlap{\kern\marginparsep\tt\@eqnlabel}%
\global\let\@eqnlabel\@vacuum}  }

\def\numberbysection{\@addtoreset{equation}{section}
        \def\theequation{\thesection.\arabic{equation}}}

\def\underline#1{\relax\ifmmode\@@underline#1\else
        $\@@underline{\hbox{#1}}$\relax\fi}

\def\titlepage{\@restonecolfalse\if@twocolumn\@restonecoltrue\onecolumn
     \else \newpage \fi \thispagestyle{empty}\c@page\z@
        \def\thefootnote{\fnsymbol{footnote}} }

\def\endtitlepage{\if@restonecol\twocolumn \else  \fi
        \def\thefootnote{\arabic{footnote}}
        \setcounter{footnote}{0}}  %\c@footnote\z@ }
\relax

\numberbysection
\hybrid

\def\beq{\begin{equation}}
\def\eeq{\end{equation}}
\def\p{\partial}

\newtheorem{th}{Theorem}[section]

\newtheorem{cor}{Corollary}[section]
\newtheorem{lem}{Lemma}[section]

\begin{document}

\begin{titlepage}

\title{Finite-gap difference operators
with elliptic coefficients and their spectral curves}

\author{A. Zabrodin
\thanks{Joint Institute of Chemical Physics, Kosygina str. 4, 117334,
Moscow, Russia and ITEP, 117259, Moscow, Russia}}
\date{October 1999}

\maketitle

\begin{abstract}

We review recent results on the finite-gap properties
of difference operators with elliptic coefficients and
give explicit characterization of spectral curves for
difference analogues of the higher Lam\'e operators.
This curve parametrizes double-\-Bloch solutions
to the difference Lam\'e equation.
The curve depends on
a positive integer number $\ell$, related to
its genus $g$ by $g=2\ell$, and two continuous
parameters: the lattice spacing $\eta$ and the modular
parameter $\tau$. Isospectral deformations of the
difference Lam\'e operator under Volterra flows
are also discussed.

\end{abstract}

%\vfill

\end{titlepage}

\section{Introduction}

The spectrum
of the Schr\"odinger operator
$-\p ^{2}_x +u(x)$ with a periodic potential
$u(x)=u(x+T)$
has a band structure: there are stable energy bands separated by
gaps. For smooth potentials, the width of gaps
rapidly decreases as energy becomes higher. However,
gaps generically occur at arbitrarily high energies,
so there are infinitely many of them.

Of particular interest are
exceptional cases, when for sufficiently high energies
there are no gaps anymore and their number is therefore finite.
Such operators are usually refered to as
algebraically integrable or finite-gap ones.
Their studies go back to classical works of the last century.
The renewed interest to the theory of finite-gap operators
is due to their role
in constructing quasi-periodic exact solutions
to non-linear integrable equations.

Among examples of the finite-gap operators, the first
and the most
familiar one is the classical Lam\'e operator
\beq
\label{Int0}
{\cal L}=-{d^2\over dx^2}+\ell (\ell +1)\wp(x +
\omega'|\omega, \omega')\,,
\eeq
where $\wp (x|\omega ,\omega')$
is the Weierstrass $\wp$-function and
$\ell$ is a parameter. The potential is a doubly-periodic
function on the complex plane with periods $2\omega$ and
$2\omega'$, where
$\mbox{Im}\,(\omega' / \omega) =\tau >0$. If $\omega$ is real
while $\omega'$ is pure imaginary, the spectral problem
is self-adjoint.
The finite gap property of
higher Lam\'e operators for integer values of
$\ell$ was established in \cite{ince}. If $\ell$ is a positive
integer, then
the Lam\'e operator has exactly
$\ell$ gaps in the spectrum.
Such a remarkable spectral property is a
signification of a hidden algebraic
symmetry which, in its turn, leads to an intimate
connection with integrable systems.

The finite-gap property becomes even more striking
for difference operators.
A natural difference analogue of the
Schr\"odinger equation has the form
\beq
\label{differ}
a(x) \Psi (x+\eta)+b(x) \Psi (x)+
c(x)\Psi (x-\eta)=E\Psi (x)\,,
\eeq
where the parameter $\eta$ is the lattice spacing.
Let us assume that $\eta$ is real, the coefficients
are real functions of $x$ and
$c(x)=a(x-\eta)$, then the problem is self-adjoint.
Let the coefficient functions be periodic with a common
period $T$: $a(x+T)=a(x)$, $b(x+T)=b(x)$.
The difference Schr\"odinger operators
with periodic coefficients exhibit much richer
spectral properties because the problem has two
competing periods ($T$ and $\eta$) rather then one.
Nevertheless, the class of finite-gap operators survives.

The structure of the spectrum of a typical
difference operator crucially depends on
whether the ratio $T/\eta$ of the two periods is rational
or irrational number. In the former
(commensurate) case one can always
set $T/\eta =Q\in {\bf Z}$ without loss of generality.
Then there are no more than $Q$ stable bands in the spectrum.
Indeed, set $\Psi (x_0 +n\eta)=\Psi_n$,
$a(x_0 +n\eta )=a_n$, etc and rewrite eq.\,(\ref{differ})
in the form
$$
a_n\Psi_{n+1} +b_n\Psi_{n} +c_n\Psi_{n-1} =E\Psi_{n,}\,,
$$
where $a_{n+Q}=a_n$, $b_{n+Q}=b_n$, $c_{n+Q}=c_n$.
Since the coefficients are periodic, one may look for
solutions in the Bloch form: $\Psi_n =e^{ik\eta n}\chi_n$,
where $\chi_n$ is $Q$-periodic and $k$ is the Bloch momentum.
Therefore, the spectral problem is reduced to the
eigenvalue problem for a hermitian
$Q$$\times$$Q$-matrix. For each real
value of the Bloch momentum $k$ the secular equation has
$Q$ real solutions $E=E_i(k)$. As $k$ sweeps over the Brillouin
zone, $E_i(k)$ sweep over the stable bands labeled by $i$.
Several neighbouring bands can merge,
so the total number of stable bands can be less
or equal to $Q$.

The latter, incommensurate case can be
practically realized as a proper limit of the former
when both numerator and denominator of the fraction
$T/\eta =Q/P$ tend to infinity. The resulting spectra can be
(and usually are) extremely complicated chaotic generations
of Cantor set type. Some of them, like those in the
Azbel-Hofstadter problem \cite{Hofst}, though of a
multifractal nature, nevertheless keep a
good deal of hidden regularity
revealed in terms of string solutions to
Bethe equations \cite{ATW}. Very little
is known on the spectra of generic type; they seem to be
completely irregular. In this paper we discuss just the opposite
case of the utmost regular spectra in the sense that
the number of bands is finite although $T/\eta$ is irrational.
Moreover, the number of bands does not really
depend on this ratio, being determined by another
(integer) parameter. The operators with this type of spectra
are true difference analogues of the finite-gap operators.

In \cite{kz}, the following
difference analogue of the Lam\'e operator (\ref{Int0})
was proposed:
\beq
L = \frac{\theta _{1}(x -\ell \eta)}
{\theta _{1}(x)}\, e^{\eta \p _x}
+\frac{\theta _{1}(x +\ell \eta)}
{\theta _{1}(x)}\, e^{-\eta \p _x}\,.
\label{Int1}
\eeq
Here $\theta_1(x)\equiv \theta_1(x|\tau)$ is the odd
Jacobi $\theta$-function, $\ell$ is a non-negative integer.
The coefficients are periodic functions with period $1$.
This operator can be made self-adjoint by
the similarity transformation $L\to g^{-1}(x)Lg(x)$
with a function $g(x)$ such that $g(x+1)=g(x)$,
so the spectrum is real.
The operator (\ref{Int1})
first appeared in a completely different
context of representations of the Sklyanin algebra
as early as in 1983 \cite{Skl}. Namely,
$L$ coincides with one
of the four generators of the Sklyanin algebra in the
functional realization found by Sklyanin.
Remarkably, the operator $L$ for positive integer values
of $\ell$ and {\it arbitrary} generic $\eta$
has $2\ell +1$ stable bands (and $2\ell$ gaps) in the spectrum.
The finite-gap property of this operator for integer
$\ell$ was proved in \cite{kz}. It was also shown \cite{kz,z}
that the Sklyanin algebra does provide a natural algebraic
framework for analyzing the spectral properties of the
operator $L$.
(A different algebraic approach to the difference
analogues of the Lam\'e operators was proposed in \cite{FV1}.)

Another similarity transformation, $\tilde L=f^{-1}Lf$,
where
\beq
\label{f}
f(x)= \prod _{j=1}^{\ell}\theta_1(x-j\eta) \,,
\eeq
makes coefficients of the difference operator
\beq
\label{tildeL}
\tilde L =e^{\eta \p_x}+
\frac{\theta_1(x+\ell \eta)\theta_1(x-(\ell +1)\eta)}
{\theta_1(x)\theta_1(x-\eta)}\,
e^{-\eta \p_x}
\eeq
double-periodic functions of $x$ with periods $1$ and $\tau$.
The limit $\eta \to 0$ gives the Lam\'e
operator (\ref{Int0}). Indeed, replacing
$x$ by $x +\frac{1}{2}\tau$ in (\ref{tildeL})
we obtain:
\beq
\label{limit}
\tilde L = 2-\eta^2({\cal L}+\mbox{const})+O(\eta^3)\,,
\;\;\;\;\;\eta \to 0\,,
\eeq
where the $\wp$-function in the ${\cal L}$ has periods $1$ and
$\tau$ (see (\ref{wei})).

Let us mention that spectral curves of the
classical Lam\'e operator (\ref{Int0}) and its
Treibich-Verdier generalizations \cite{TV} for small
values of $\ell$ were studied in \cite{En},\,\cite{Smir}.
A detailed analysis of solutions to
the difference Lam\'e equation
at $\ell =1$ was recently carried out in \cite{Ruijs}.

The paper is organized as follows. Sec.\,2 is
a continuation of the Introduction. To draw the
problem in a broader context, we discuss the general
notion of the finite-gap operator. In Sec.\,3, a
family of Bloch eigenfunctions of the operator (\ref{tildeL})
is constructed. These eigenfunctions are parametrized
by points of the spectral curve. Sec.\,4 contains
equations for the edges of bands and some examples.
In Sec.\,5 we work out an explicit relation between the
Bloch multipliers. The form of the result suggests that
some hypothetical combinatorial identities for "elliptic numbers"
may be relevant. At last, Sec.\,6 contains some remarks
on the isospectral deformations of the difference Lam\'e
operator. In this case the coefficient in (\ref{tildeL})
has more poles. The location of the poles,
however, is not arbitrary: they are constrained by locus
equations.

\section{A general view of finite-gap operators}

The key idea of the modern approach to spectra
of differential or difference operators
is to regard the solutions $\Psi (x,E)$
to the spectral problem (say, (\ref{differ})) as functions
of $E$ for any complex values of $E$ and to study
their analytic properties in $E$. In so doing
it is not necessary to assume that the problem is
self-adjoint, so the parameter $\eta$ and the
coefficients may be complex numbers.

In practice, one may try to construct a family of
eigenfunctions $\Psi \!=\! \Psi (x, E, p_1, p_2, \ldots )$
depending, apart from $E$, on a finite number of additional
parameters $p_i$. For instance, one of these could be
the Bloch momentum $k$: $\Psi (x)=e^{ikx}\chi (x,E, k)$,
where $\chi (x+T,E,k)=\chi (x,E,k)$.
Suppose such a family does exist. Then the
spectral parameters appear to be constrained by some relations
$F_{i}(E,p_1, p_2, \ldots )=0$, so that only one of the
parameters is independent. These relations define
a complex curve (a Riemann surface)
in the parameter space called
{\it the spectral curve}. The true spectral
parameter is a point of the curve.
This is the proper mathematical
formulation of the dispersion law $E=E(k)$.
Usually, this function is multi-valued. It becomes
single-valued on the spectral curve (when the latter is
well-defined). Moreover, the solution $\Psi (x,E)$ of the spectral
problem also becomes a single-valued function on the
spectral curve.

In the case of the second-order difference operators
the spectral problem (\ref{differ}) has no more than two
linearly independent solutions. In other words, the function
$E$ can take any of its values at most twice. The existence
of such a function implies
that the spectral curve is a hyperelliptic curve. Any
hyperelliptic curve can be represented in the form
a two-sheet covering of the complex plane of the variable $E$:
\beq
\label{hyp}
w^2=\prod_{i}(E-E_i)\,,
\eeq
where $E_i$ are called {\it branch points}.
For $E\!=\!E_i$ the equation (\ref{differ}) has only one linearly
independent solution.
The curve (\ref{hyp}) is well-defined
if the set of branch points is finite.
Then the curve is algebraic and
has finite genus. Equivalently, this means that
there exists a difference operator $W$ such that $W$
can not be represented as a polynomial function
of the difference operator in the right hand side of
(\ref{differ}), and that commutes with this operator.
In this case
they have a set of common eigenfunctions.
Equation (\ref{hyp}) is then lifted to the operator
relation. The parameter $w$ is
the eigenvalue of the operator $W$ on the common
eigenfunction
$\Psi (x,E)$: $W\Psi (x,E)=w\Psi (x,E)$.

For self-adjoint spectral problems the branch points $E_i$
are real numbers $E_1 < E_2 < E_3 <\ldots $.
In the stable bands the Bloch momentum takes real values.
The stable bands are segments of the
real line $[E_{2i+1} ,\, E_{2i+2}]$, $i=0,1, \ldots$,
so the branch points are just edges of bands.
At the edges of bands, the Bloch solutions
$\Psi (x, E_i)$ are periodic or anti-periodic.

The notion of the spectral curve
is really useful in the exceptional case of the spectra
of regular type mentioned in the Introduction.
Recall that
in the typical case the set of branch points may be
even uncountable, so the spectral curve, defined as above,
does not have sense.
The very existence of a well-defined spectral curve
of finite genus is the precise
characterization of the finite-gap operators.

Then it is natural to address the inverse problem:
given a hyperelliptic curve of finite genus
regarded as a spectral curve of some difference
operator, to find coefficients of this operator,
i.e., the functions $a(x)$, $b(x)$ in (\ref{differ}).
In this way, one is able to construct a representative family
of finite-gap operators \cite{fingap1,fingap2}. For difference
operators this was done in \cite{dattan,mum,kr1}.
The coefficients of the operators can be expressed
through Riemann's theta-functions
associated with the curve.

The curve itself does not determine the operator uniquely.
There is a remaining finite-parametric freedom
which can be fixed by some additional data on the curve
(essentially, a number of marked points).
In other words, any finite-gap operator admits a class
of isospectral deformations.
The coefficients of the operator, with respect to the
isospectral flows, obey certain non-linear integrable
equations.

Looking for formulas more effective than a bunch of Riemann's
theta-fucntions,
one may inquire whether they can be expressed
through simpler functions, for instance elliptic ones.
For a particular class of curves, namely, for special
coverings of elliptic curves,
this is indeed possible (see
e.g. \cite{book}).
The Riemann theta-function associated with
such a curve factorizes
into a product of Jacobi $\theta$-functions, so the
coefficients of the operator become elliptic functions.
Moreover, all the family of isospectral deformations of the
operator enjoys the same property.

An important example of this phenomenon in the
differential set-up is provided by
the Lam\'e operator (\ref{Int0})
for $\ell \in {\bf Z}_+$ and its isospectral deformations
$-\p _{x}^{2}+u(x)$ with
\beq
\label{iso1}
u(x)=
2 \sum_{j=1}^{\ell (\ell \!+\!1)/2}\wp(x -x_j )
+\,\,\mbox{const}\,.
%\;\;\;\;\;\; N=\frac{1}{2}\ell (\ell +1)
\eeq
The Lam\'e operator itself corresponds
to a very degenerate configuration when
all the poles sit in one and the same point.
The isospectral flows are the
flows of the KdV hierarchy \cite{fingap1} for the potential
$u(x)$. Solving say the KdV equation
$\dot u=6uu' -u'''$ for $u=u(x,t)$ with the initial condition
$u(x,0)=\ell (\ell \!+\!1)\wp (x-x_0)$, we get a family
of Schr\"odinger operators with elliptic potential which have
the same spectral curve as the Lam\'e operator.
The poles $x_j$ (and the constant term)
in (\ref{iso1}) become $t$-dependent.
By a direct substitution to the KdV equation,
it can be shown \cite{elkdv} that they
are constrained by the conditions
\beq
\label{loc1}
\sum _{j=1, \neq i}^{\ell (\ell \!+\!1)/2} \wp '(x_i -x_j ) =0\,,
\;\;\;\;\;\;i=1,2, \ldots , \ell (\ell \!+\!1)/2\,,
\eeq
and obey the differential equations
\beq
\label{iso2}
\dot x_j =-12
\sum _{k=1, \neq j}^{\ell (\ell \!+\!1)/2} \wp (x_j -x_k )\,.
\eeq
Eqs.\,(\ref{loc1}) are the famous equations defining the equilibrium
locus of the elliptic Calogero-Moser system of particles.
From the general theory which connects
the pole dynamics of elliptic solutions of non-linear integrable
equations with systems of Calogero-Moser
type \cite{elkdv,krich} it follows that
the connected component of the locus is parametrized by the Jacobian of
the spectral curve of the Lam\'e operator.
So, it is an $\ell$-dimensional submanifold, spanned
by higher Calogero-Moser flows, in the
$\frac{1}{2}\ell (\ell \!+\!1)$-dimensional configuration
space with coordinates $x_j$.

In Sec.\,6, we present analogues of equations
(\ref{iso1}), (\ref{loc1}) and (\ref{iso2})
in the difference set-up. The isospectral flows
are connected with eliptic solutions to the Volterra
hierarchy.

\section{Double-Bloch eigenfunctions
of the difference Lam\'e operator and the spectral curve}

In this section we study
Bloch eigenfunctions of the difference Lam\'e
operator. Following the general scheme outlined at the beginning
of Sec.\,2, we construct a family of eigenfunctions depending
on $E$ and two additional spectral parameters. All the
three parameters are constrained by two equations which
define the spectral curve.

Consider the eigenvalue equation for the operator $\tilde L$
(\ref{tildeL}):
\beq
\label{db2}
\psi(x+\eta)+
\frac{\theta_1(x+\ell \eta)\theta_1(x-(\ell +1)\eta)}
{\theta_1(x)\theta_1(x-\eta)}\,
\psi(x-\eta)=E\psi(x)\,.
\eeq
The coefficient function is
double-periodic. Therefore,
it is natural to look for solutions in the class of
{\it double-Bloch functions} \cite{kz}, i.e., such that
$\psi (x+1)=B_1 \psi(x)$,
$\psi (x+\tau )=B_{\tau} \psi(x)$ with some
constants $B_1, B_{\tau}$. These are going to be the additional
parameters $p_i$ from the general scheme of Sec.\,2.

Consider the function
\beq
\label{Phi}
\Phi (x,\zeta)=
\frac{\theta _{1}(\zeta +x)}{\theta _{1}(x)\theta _{1}(\zeta)}\,.
\eeq
Its monodromy properties in $x$ are
$\Phi(x+1, \zeta)=\Phi (x,\zeta)$,
$\Phi(x+\tau, \zeta)=e^{-2\pi i \zeta}\Phi (x,\zeta)$,
i.e., it is a double-Bloch function. Moreover, it is the simplest
non-trivial (i.e., different from the exponential function)
double-Bloch function since it has only one pole. This
function serves as a building block for more general
double-Bloch functions.

Let $\ell$ be a positive integer. We employ the following
double-Bloch ansatz for the $\psi$:
\beq
\label{db3}
\psi(x)=K^{x/\eta}\sum_{j=1}^{\ell}
s_j(\zeta, K, E) \Phi (x -j\eta, \zeta)\,,
\eeq
where $\zeta ,K$ parametrize the Bloch multipliers of the
function $\psi (x)$:
$B_1=K^{\frac{1}{\eta}}$, $B_{\tau}=
K^{\frac{\tau}{\eta}}e^{-2\pi i \zeta}$.
The coefficients $s_j$ depend on the indicated parameters
only.

Substituting (\ref{db3}) into
(\ref{db2}) and computing the residues at the points
$x=j\eta,\ j=0,\ldots,\ell$, we get $\ell +1$ linear equations
\beq
\sum_{j=1}^{\ell} M_{ij} s_j=0,
\;\;\;\;\;\; i=0,1,\ldots, \ell \,,
\label{db4}
\eeq
for  $\ell$ unknowns $s_j$.
Matrix elements $M_{ij}$ of this system are given by the
formula

\begin{eqnarray}
\label{db5}
M_{ij}&=&K\delta_{i,j\!-\!1}-E\delta_{i,j}+K^{-1}
\frac{\theta_1((j+\ell +1)\eta)
\theta_1((j-\ell )\eta)}{\theta_1((j+1)\eta)\theta_1(j\eta)}
\delta_{i,j\!+\!1}+
\nonumber\\
&+&K^{-1}\frac{\theta_1(\zeta \!-\!(j\!-\!i\!+\!1)
\eta)}{\theta_1(\zeta)}\,
\frac{\theta_1((i+\ell)\eta)
\theta_1((i\!-\!\ell \!-\!1)\eta)}{\theta_1(\eta)
\theta_1((j\!-\!i\!+\!1)\eta)}
(\delta_{i,0}-\delta_{i,1})\,.
\end{eqnarray}
Here $i=0,1, \ldots , \ell$,
$j=1,2, \ldots , \ell$.
The overdetermined system (\ref{db4}) has nontrivial solutions if
and only if
rank of the rectangular matrix $M_{ij}$ is less
than $\ell$. Let $M^{(0)}$ and $M^{(1)}$ be
$\ell \times  \ell$ matrices obtained from $M$ by deleting the rows
with $i=0$ and $i=1$,
respectively. Then the values of three parameters
$\zeta,K,E$ for which eq.\,(\ref{db4})
has solutions of the form (\ref{db3})
are determined by the system of two equations:
$\det M^{(0)}=\det M^{(1)}=0$. They indeed define
a curve.

To obtain the explicit form of these equations,
we expand the determinants
with respect to the first row.
This yields an explicit characterization
of the spectral curve summarized in the theorem below.
Hereafter, the
"elliptic factorial" and "elliptic binomial"
notation is convenient:
\beq
\label{binom}
\begin{array}{l}
\displaystyle{[n]!=\prod_{j=1}^{n}[j]}\,,
\;\;\;\;\;\;[j]\equiv \theta_1(j\eta)/\theta_1(\eta)\,,
\\ \\
\left [ \begin{array}{c}n\\m\end{array}\right ]
\equiv \displaystyle{\frac{[n]!}{[m]![n-m]!}}\,.
\end{array}
\eeq

\begin{th}
The difference Lam\'e equation (\ref{db2})
has double-Bloch solutions of the form (\ref{db3})
if and only if the spectral parameters $\zeta, K, E$
obey the equations
\beq
\begin{array}{l}
\displaystyle{\sum_{j=0}^{\ell}(-1)^j K^{-j}
\theta_1(\zeta-j\eta)\left [
\begin{array}{c}\ell\\ j \end{array}\right ]}
A^{(\ell )}_{j}(E)=0\,, \\ \\
\displaystyle{\sum_{j=0}^{\ell +1}(-1)^j K^{-j}
\theta_1(\zeta-j\eta)[j-1] \left [
\begin{array}{c}\ell \!+\!1\\ j \end{array}\right ]}
A^{(\ell )}_{|j-1|}(E)=0\,,
\end{array}
\label{db7}
\eeq
where $A^{(\ell )}_{j}(E)$ are polynomials of $(\ell -j)$-th degree
explicitly given by the determinant formula
\beq
\label{AE}
A^{(\ell )}_{\ell -s}(E)=\left [\begin{array}{c}\ell \\ s
\end{array}\right ]
\left [\begin{array}{c}2\ell \\ s
\end{array}\right ]^{-1} \det \left (
E\delta _{i,j}+
\frac{[-i]}{[\ell \!+\!1\!-\!i]}\delta _{i,j-1}
+\frac{[2\ell \!+\!2\!-\!i]}{[\ell \!+\!1\!-\!i]}\delta _{i,j+1}
\right )_{1\leq i,j\leq s}
\eeq
(here $0\leq s\leq \ell$), $A^{(\ell )}_{\ell}=1$.
\end{th}

Let us list some useful properties of the polynomials
$A^{(\ell )}_{j}(E)$.
First, they obey the recurrence
relation:
\beq
\label{bax2}
A^{(\ell )}_{\ell \!-\!s\!-\!1}(E)=
\frac{[\ell \!-\!s]}{[2\ell \!-\!s]}EA^{(\ell )}_{\ell \!-\!s }(E)
+\frac{[s]}{[2\ell \!-\!s]}A^{(\ell )}_{\ell \!-\!s\!+\!1}(E)
\eeq
with the initial condition $A^{(\ell )}_{\ell }(E)=1$,
$A^{(\ell )}_{\ell -1}(E)=([\ell ]/[2\ell ])E$.
Next, it is clear from (\ref{bax2}) that
\beq
\label{parity}
A^{(\ell )}_{\ell -s}(-E)=
(-1)^s A^{(\ell )}_{\ell -s}(E)\,,
\;\;\;\;\;\;
0\leq s\leq \ell\,.
\eeq

The equations (\ref{db7})
define a Riemann surface $\tilde \Gamma$, which covers the
complex plane. The monodromy properties of the $\theta$-function
(see (\ref{periods}) in the Appendix) make it clear that
this surface is invariant under the transformation
\beq
\label{db8}
\zeta \longmapsto  \zeta +\tau\,,
\;\;\;\;\;\;
K \longmapsto Ke^{2\pi i\eta}\,.
\eeq
The factor of the $\tilde \Gamma$ over this transformation
is an algebraic curve $\Gamma$, which is
a ramified covering of the elliptic
curve with periods $1$, $\tau$. It is clear
from (\ref{db7}), (\ref{parity}) that the curve
admits the involution
\beq
\label{db81}
(\zeta , K,E)\longmapsto (\zeta , -K,-E)\,,
\eeq
so the spectrum is symmetric with respect
to the reflection $E\rightarrow -E$.
Another result of \cite{kz}, which is not so easy to see
from (\ref{db7}), is that the curve $\Gamma$ is at the same
time a hyperelliptic curve.
\begin{th}
The curve $\Gamma$
is a hyperelliptic curve of genus $g=2\ell$. The hyperelliptic
involution is given by
\beq
\label{db9}
(\zeta , \,K, \,E) \; \longmapsto \;
(2N\eta -\zeta, \,K^{-1}, E)\,,
\;\;\;\;\;\;N=\frac{1}{2}\ell (\ell +1)\,.
\eeq
The points $P=(\zeta , K, E) \in \Gamma$ of the curve
parametrize double-Bloch solutions
$\psi (x)=\psi (x,P)$ to eq.\,(\ref{db2}), and the
solution $\psi (x,P)$
corresponding to each point
$P\in \Gamma$ is unique up to a constant multiplier.
\end{th}
For the proof see \cite{kz}. Here we give a few remarks.
The involution (\ref{db9}) looks best in terms of
the function $\Psi (x)=
\psi (x)\prod_{j=1}^{\ell}\theta_1(x-j\eta)$ which satisfies
the eigenvalue equation $L\Psi =E\Psi$ with the $L$ as in
(\ref{Int1}) (cf. (\ref{f})). Then the hyperelliptic
involution simply takes $\Psi (x)$ to $\Psi (-x)$.
The genus of the
curve can be found by counting the number of the fixed points
of this involution. At the fixed points the two solutions,
$\Psi (x)$ and $\Psi (-x)$, are linearly dependent:
$\Psi (-x)=r\Psi (x)$. Writing out the eigenvalue equation
at $x=0$, we obtain the necessary
condition $\Psi (-\eta)=\Psi (\eta)$, so $r=1$. The ansatz
(\ref{db3}) for $\psi$ is equivalent to the ansatz
$$
\Psi (x)=K^{x/\eta}\prod_{j=1}^{\ell}\theta_1 (x+y_j)
$$
with $\sum_j y_j =\zeta$. At $K=1$,
the dimension of the linear space of
even functions of this form is known to be
$2\ell +1$. Adding images of the fixed points under the
involution (\ref{db81}), we eventually get $4\ell +2$ fixed
points, so, by the Riemann-Hurwitz formula,
the genus is equal to $2\ell$.

Taking into account the symmetry $E\to -E$,
we can write the equation of the hyperelliptic
curve in the standard form:
\beq
\label{15a}
w^2=\prod_{i=1}^{2\ell +1}(E^2-E_{i}^{2})\,.
\eeq
The hyperelliptic involution takes $(w,E)$ to $(-w,E)$.

In (\ref{15a}), $w$ is the eigenvalue of a non-trivial operator
$W$ commuting with $L$ on their common eigenfunction.
The explicit form of the operator $W$ was found in \cite{FV2}:
$$
W = \varphi _{\ell}(x)\sum_{k=0}^{2\ell +1} (-1)^k
\left [ \begin{array}{c} 2\ell +1 \\ k \end{array} \right ]
\frac{ \theta _{1}(x+
(2\ell -2k +1)\eta )}{\prod_{j=0}^{2\ell -k+1}
\theta_1(x+j\eta )
\prod_{j'=1}^{k}\theta_1(x-j'\eta )}
e^{(2\ell -2k +1)\eta \p_x }\,,
$$
where $\varphi _{\ell}(x)=\prod _{j=0}^{2\ell }
\theta_1(x+(j-\ell )\eta )$.

Let us conclude this section by examining the behaviour
of the spectral curve in the vicinity of its "infinite points",
i.e., the points at which the function $E$ has poles.
From (\ref{db7}) we
conclude that there are two such points:
$\infty_{+}=(\zeta \to 0,\, K\to \infty ,\, E\to \infty )$
and $\infty_{-}=
(\zeta \to 2N\eta ,\, K\to 0 ,\, E\to \infty )$.
In the neighbourhood of $\infty _{\pm}$
$E=K^{\pm 1} +o(K^{\pm 1})$. In terms of the variables $(w,E)$
these points are $\infty _{\pm}=(w\to \pm \infty ,\, E\to \infty )$.

\section{Edges of bands}

The edges of bands $\pm E_i$, i.e., the
branch points of the two-sheet covering (\ref{15a}),
are values of the function $E=E(P)$ at the fixed points of the
hyperelliptic involution.
As is clear from (\ref{db9}), the fixed points lie
above the points $\zeta = N\eta +\omega_a$, where $\omega_a$
are the half-periods: $\omega_1=0$, $\omega_2 =\frac{1}{2}$,
$\omega_3 =\frac{1}{2}(1+\tau )$,
$\omega_4 =\frac{1}{2}\tau $. The corresponding values of $K$
are determined from (\ref{db8}). Then the set of the
branch points $E_i$ is fixed by Theorem 3.1.

\begin{cor}
Let ${\cal E}_a$, $a=1,\ldots , 4$ be the set of common
roots of the polynomial equations
\beq
\begin{array}{l}
\displaystyle{\sum_{j=0}^{\ell}
(-1)^j \theta_a((N-j)\eta)\left [
\begin{array}{c}\ell\\ j \end{array}\right ]}
A^{(\ell )}_{j}(E)=0\,, \\ \\
\displaystyle{\sum_{j=0}^{\ell +1}
(-1)^j \theta_a((N-j)\eta)[j-1] \left [
\begin{array}{c}\ell \!+\!1\\ j \end{array}\right ]}
A^{(\ell )}_{|j-1|}(E)=0\,,
\end{array}
\label{db10}
\eeq
where
$\theta_a$ are Jacobi $\theta$-functions\footnote{Definitions
and transformation properties
of the Jacobi $\theta$-functions
$\theta_a(x|\tau )$, $a=1,2,3,4$, are listed in the Appendix.
For brevity, we write $\theta _a(x|\tau )
\equiv \theta _a(x)$.}.
Then the set of the edges of bands $\pm E_i$
is the union of $\,\bigcup_{a=1}^{4}{\cal E}_a$ and its
image under the reflection $E\rightarrow -E$.
\end{cor}

Let us give two examples. At $\ell =1$, the set ${\cal E}_1$
is empty while ${\cal E}_{a}$ for $a=2,3,4$ contains one point.
From (\ref{db10}) we find:
$$
E_{\alpha}=2\,
\frac{\theta _{\beta +1}(\eta )
\theta _{\gamma +1}(\eta )}{\theta _{\beta +1}(0)
\theta _{\gamma +1}(0)}\,,
$$
where $\{\alpha , \,\beta ,\, \gamma \}$ is any cyclic
permutation of $\{ 1,\, 2, \, 3 \}$.
At $\ell =2$, the set ${\cal E}_1$ has two elements $E_1 , E_2$
obtained as solutions of the quadratic equation
$[2]E^2 +[2]^3 E+2[4]=0$,
so that
$$
E_{1,2}=\frac{1}{2}\left (
\frac{\theta_{1}^{2}(2\eta )}{\theta_{1}^{2}(\eta )}\pm
\sqrt{
\frac{\theta_{1}^{4}(2\eta )}{\theta_{1}^{4}(\eta )}
-8 \frac{\theta_{1}(4\eta )}{\theta_{1}(2\eta )} } \,\right )\,.
$$
For $a=2,3,4$ each set
${\cal E}_a$ has one element $E_{a+1}$:
$$
E_{a+1}=
\frac{\theta _{1}(2\eta )
\theta _{a}(2\eta )}{\theta _{1}(\eta )
\theta _{a}(\eta )}\,.
$$

In general, it is possible to prove \cite{kz} that
\beq
\label{E1234} \#({\cal E}_{1})=
\left \{
\begin{array}{lll}
\frac{1}{2}(\ell \! -\! 1)\,,&
\!\!& \ell \;\; \mbox{odd}
\\ \\
\frac{1}{2}\ell \!+\!1\,, &
\!\! &\ell \;\; \mbox{even}
\end{array}\right.
\;\;\;\;\; \#({\cal E}_{2,3,4})=
\left \{
\begin{array}{lll}
\frac{1}{2}(\ell \!+\!1)\,, &
\!\! &\ell \;\; \mbox{odd}
\\ \\
\frac{1}{2}\, \ell \,, &
\!\! &\ell \;\; \mbox{even}\,.
\end{array}\right.
\eeq
Note that $\# \Bigl (\cup _{a=1}^{4}{\cal E}_{a}\Bigr ) =2\ell +1$
that agrees with (\ref{15a}).

\section{Relation between the Bloch multipliers}

To simplify
equations of the curve (\ref{db7}), one can try to eliminate
one of the variables and obtain a single equation
for the other two. Here we show how to eliminate $E$.
This leads to a closed relation between
the two Bloch multipliers of the function
(\ref{db3}) (parametrized through $\zeta$ and $K$).
Its form (see (\ref{eq4}),
(\ref{Cj}) below) suggests an interpretation
in terms of hypothetical
combinatorial identities for elliptic numbers.

At the first glance, the elimination of $E$ from eqs.\,(\ref{db7})
is hardly possible at all.
Nevertheless, there is an alternative argument leading directly
to the relation between the Bloch multipliers.
Here it is more convenient to deal with the difference
Lam\'e operator in the gauge equivalent form (\ref{Int1}).
Our construction is based on the following
simple lemma \cite{z}.
\begin{lem}
Let $\Psi (x)$ be any solution to the equation
\beq
\label{eq1}
\frac{\theta _{1}(x -\ell \eta)}{\theta _{1}(x)}
\Psi (x+\eta ) +
\frac{\theta _{1}(x +\ell \eta)}{\theta _{1}(x)}
\Psi (x-\eta )
=E\Psi (x)
\eeq
in the class of entire functions on the complex plane of
the variable $x$, then
\beq
\label{cond}
\Psi (j\eta )=\Psi (-j\eta )\,,
\;\;\;\;\;\;\;\;\; j =1,2, \ldots , \,\ell\,.
\eeq
\end{lem}
This assertion follows from the specific location
of zeros and poles of the coefficients
of eq.\,(\ref{eq1}). Indeed, putting $x=0$ in (\ref{eq1}),
we have $\Psi (\eta)=\Psi (-\eta)$. The proof can be
completed by induction. At $x=\pm \ell \eta$ one of the coefficients
in the l.h.s. of (\ref{eq1}) vanishes, so
the chain of relations (\ref{cond}) truncates at $j=\ell$.

Remarkably, the conditions (\ref{cond})
and the ansatz
\beq
\label{psi1}
\Psi (x)=K^{x/\eta}
\left (\prod _{j=1}^{\ell}\theta _{1}(x-j\eta )\right )
\sum _{m=1}^{\ell}s_{m}(K,\zeta )\Phi (x -m\eta , \, \zeta )
\eeq
for $\Psi$ (equivalent to the ansatz (\ref{db3})
for $\psi$) with the same function $\Phi (x,z)$ given by
(\ref{Phi}) allow one
to find the relation between the
Bloch multipliers even without explicit use of
the difference Lam\'e equation (\ref{eq1}).
Plugging (\ref{psi1}) into (\ref{cond}), we obtain $\ell$
equalities (for $m=1,2, \ldots , \ell$):
\beq
\label{eq2}
K^m s_m =(-1)^{\ell}K^{-m}\theta _{1}(2m\eta )
\left ( \prod _{j=1, \neq m}^{\ell}
\frac{\theta _{1}((m+j)\eta )}
{\theta _{1}((m-j)\eta )}\right )
\sum _{n=1}^{\ell} \Phi \Bigl (-(m+n)\eta , \, \zeta \Bigr )\, s_n\,.
\eeq
This is a system of linear homogeneous equations for $s_n$.
It has nontrivial solutions if and only if its determinant
is equal to zero, whence we obtain
the equation connecting $\zeta$ and $K$:
\beq
\label{eq3}
\det \Bigl ( K^{2m}\delta _{mn}+G_{mn}(\zeta )
\Bigr )_{1\leq m,n\leq \ell}=0\,,
\eeq
where
$$
G_{mn}(\zeta )=
(-1)^{\ell +1}[2m]
\left ( \prod _{j=1, \neq m}^{\ell}
\frac{[m+j]}{[m-j]}\right )
\Phi \bigl (-(m+n)\eta , \, \zeta \bigr )\,.
$$
This equation defines a curve, which is
the image of the spectral curve $\Gamma$ under the
projection that takes
$(\zeta , K, E)$ to $(\zeta , K)$.

The equation of the spectral curve (\ref{eq3})
can be represented in the form
\beq
\label{eq4}
\sum _{j=0}^{N}(-1)^j C^{(\ell)}_{j}
(\eta )\theta _{1}(\zeta -2j\eta)
K^{2(N-j)}=0\,,
\eeq
where $N=\frac{1}{2}\ell (\ell +1)$ and $C_{j}^{(\ell)}(\eta)$ are
some coefficients depending only on $\eta$ and $\tau$ such that
$C_{j}^{(\ell)}(\eta)=C_{N-j}^{(\ell)}(\eta)$,
$C_{0}^{(\ell)}(\eta)=1$.
To see this, we
expand the determinant (\ref{eq3}) in powers of $K$
with the help of the identity
$$
\det \left (\frac{\theta_{1}
(x_i+x_j+\zeta)}{\theta_{1}(x_i+x_j)}
\right )_{1\leq i,j\leq n}\!\!=\,
\frac{\theta _{1}^{n-1}(\zeta)\theta _{1}(\zeta
+2\sum_{i=1}^{n}x_i)}{\prod_{i=1}^{n}\theta_{1}(2x_i)}\,
\prod_{i<j}^{n}
\frac{\theta_{1}^{2}(x_i-x_j)}{\theta_{1}^{2}(x_i+x_j)}
$$
(a particular case of the formula for the elliptic
Cauchy determinant). Let $\Lambda$
be the set $\{1,2,\ldots ,\ell\}$.
For any subset
$J\subseteq \Lambda$,
$\Lambda \setminus J$ is its complement, and
$\sigma (J)  =\sum_{j\in J}j$.
Evaluating the elliptic Cauchy determinant at
$x_n=-2\eta n$, we get:
\begin{eqnarray}
\label{e1}
&&\det \Bigl ( K^{2m}\delta _{mn}+G_{mn}(\zeta )
\Bigr )_{1\leq m,n\leq \ell}
\nonumber\\
&=& \sum _{J\subseteq \Lambda}
\frac{\theta _{1}(\zeta -2\sigma (J)\eta)}{\theta _{1}(\zeta)}
K^{2N-2\sigma (J)}(-1)^{\sigma (J)}
\prod_{k\in J}\prod_{k'\in \Lambda \setminus J}
\frac{[k+k']}{[|k-k'|]}\,.
\end{eqnarray}
Thus, the coefficient $C_{j}^{(\ell)}$ reads
\beq
\label{Cj}
C_{j}^{(\ell)}=
\sum_{J\subseteq \Lambda , \sigma (J)=j}
\prod_{k\in J}\prod_{k'\in \Lambda \setminus J}
\frac{[k+k']}{[|k-k'|]}\,.
\eeq
The symmetry $j\leftrightarrow N-j$ is now transparent.
Note that the sum in (\ref{Cj}) runs
over partitions of the number $j$ into {\it distinct} parts not
exceeding $\ell$.

We note that for any elliptic module $\tau$
\beq
\label{bin}
\lim _{\eta \to 0}C_{j}^{(\ell )}=
\left ( \begin{array}{c}N\\j \end{array}\right )
\eeq
(the usual binomial coefficient). To see this, consider
the limiting case $\tau \to i\infty$. We have:
$\exp (-\pi i \tau /4)\theta_1(x|\tau)\to 2\sin (\pi x)$
as $\tau \to i\infty$ (see (\ref{infp})),
so $$[j]\longrightarrow
\frac{\sin (\pi \eta j)}{\sin (\pi \eta )}\equiv
(j)_q\,,
\;\;\;\;\;\;\; q=e^{2\pi i \eta }\,.
$$
Then (\ref{bin}) follows from the identity \cite{DK}
\beq
\label{schur}
\sum _{J\subseteq \Lambda}z^{\sigma (J)}
\prod_{k\in J}\prod_{k'\in \Lambda \setminus J}
\frac{(k+k')_q}{(|k-k'|)_q}=
\prod_{1\leq j\leq k\leq \ell}(1+zq^{j+k-\ell -1})
\eeq
which is a specialization of the $C_{\ell}$-type Weyl
denominator formula\footnote{I am grateful to A.N.Kirillov
for pointing out this identity and drawing my attention
to the paper \cite{DK}.}. A detailed combinatorial
analysis of the
limiting cases $\tau \to i\infty$ or $\tau \to 0$
of the difference Lam\'e operators can be found in \cite{DK}.

\section{Isospectral deformations of the difference
Lam\'e operator and locus equations}

Finally, let us comment on isospectral deformations
of the difference Lam\'e operator. We are going to
present difference analogues of the operator (\ref{iso1})
and of the locus equations (\ref{loc1}).

Instead of (\ref{db2}), consider the equation
\beq
\label{more1}
\psi (x+\eta )+c(x)\psi (x-\eta ) =E\psi (x)
\eeq
with a more general coefficient $c(x)$, which is an elliptic
function represented in the form
\beq
\label{more2}
c(x)=\frac{\rho (x+\eta )\rho (x-2\eta )}{\rho (x)\rho (x-\eta )}\,,
\;\;\;\;\;
\rho (x)=\prod _{j=1}^{\ell (\ell \!+\!1)/2}
\theta _{1}(x-x_j )\,.
\eeq
Note that in the case of the difference Lam\'e operator
the configuration of zeros of the $\rho (x)$ is very specific:
\beq
\label{rho0}
\rho (x)=\prod _{1\leq j\leq k \leq \ell }
\theta _{1}(x +(j+k-\ell -1)\eta )\,,
\eeq
so all but two cancel in the $c(x)$.

The isospectral flows are the flows of the Volterra
hierarchy for the $c(x)$.
The first equation of the hierarchy,
\beq
\label{volt}
\frac{\p c(x)}{\p t}=-c(x)(c(x+\eta )-c(x-\eta ))\,,
\eeq
is the compatibility condition
of the spectral problem (\ref{more1}) and the linear problem
$$
\frac{\p \psi (x)}{\p t}=c(x)c(x-\eta )\psi (x-2\eta )\,.
$$
Recall that changing the variables as
$t\to \frac{1}{3}\eta^3 t$, $x\to x-2\eta t$
and setting $c(x)=1-\eta^2 u(x)$, one gets the
KdV equation for $u$ as $\eta \to 0$.

Substituting the pole ansatz (\ref{more2}) into the
Volterra equation and requiring the residues at the poles
to be zero, we get the following two systems of
equations ($j=1,2,\ldots , \frac{1}{2}\ell (\ell \!+\!1)$):
\beq
\label{more3}
\begin{array}{l}
\displaystyle{
\frac{\theta_{1}'(0)}{\theta _{1}(2\eta )}\,\dot x_j=
\prod _{k=1, \neq j}^{\ell (\ell \!+\!1)/2}
\frac{\theta_1(x_j -x_k +2\eta )
\theta_1(x_j -x_k -\eta )}{\theta_1(x_j -x_k +\eta )
\theta_1(x_j -x_k)}  }\,,
\\ \\
\displaystyle{
\frac{\theta_{1}'(0)}{\theta _{1}(2\eta )}\,\dot x_j=
\prod _{k=1, \neq j}^{\ell (\ell \!+\!1)/2}
\frac{\theta_1(x_j -x_k -2\eta )
\theta_1(x_j -x_k +\eta )}{\theta_1(x_j -x_k -\eta )
\theta_1(x_j -x_k)}  }\,,
\end{array}
\eeq
where $\dot x_j=\p _{t}x_j$.
Solving these equations, or,
equivalently, the Volterra equation with the initial condition
(\ref{rho0}), one arrives at a family of isospectral deformations
of the difference Lam\'e operator.
The two systems (\ref{more3}) must be
satisfied simultaneously, therefore, the right hand sides
are identical.
Whence we obtain the necessary
conditions for solutions to exist:
\beq
\label{loc3}
\prod _{k=1, \neq j}^{\ell (\ell \!+\!1)/2}
\frac{\theta_1(x_j -x_k +2\eta )
\theta_1^2(x_j -x_k -\eta )}{\theta_1(x_j -x_k -2\eta )
\theta_1^2(x_j -x_k +\eta )}  =1
\eeq
that is the difference analogues of (\ref{loc1}). The results
of \cite{kz} imply that these equations define an equilibrium
locus of the Ruijsenaars-Schneider model.
Like in the differential case, the locus is not compact.
Its closure contains, in particular, the point corresponding
to the degenerate configuration (\ref{rho0}).
From the general
arguments of ref.\,\cite{kz} it follows that the connected
component of the locus that contains
this point at the boundary is $\ell$-dimensional.
Expanding (\ref{loc3}) in $\eta \to 0$, we do get (\ref{loc1}).

\section*{Acknowledgments}
The author would like to thank Professors M.Kashiwara and
T.Miwa for the opportunity to present these results
on the workshop "Physical Combinatorics".
Discussions with I.M.Krichever,
A.D.Mi\-ro\-nov, T.Ta\-ke\-be and P.B.Wieg\-mann
are gratefully acknowledged.
This work was supported in part by RFBR grant 98-01-00344
and by grant for support of scientific schools.

\section*{Appendix. Theta-functions}
%\addcontentsline{toc}{section}{Appendix}
\def\theequation{A\arabic{equation}}
\setcounter{equation}{0}

We use the following definition of the
Jacobi $\theta$-functions:
\beq
\begin{array}{l}
\theta _1(x|\tau)=-\displaystyle{\sum _{k\in {\bf Z}}}
\exp \left (
\pi i \tau (k+\frac{1}{2})^2 +2\pi i
(x+\frac{1}{2})(k+\frac{1}{2})\right ),
\\ \\
\theta _2(x|\tau)=\displaystyle{\sum _{k\in {\bf Z}}}
\exp \left (
\pi i \tau (k+\frac{1}{2})^2 +2\pi i
x(k+\frac{1}{2})\right ),
\\ \\
\theta _3(x|\tau)=\displaystyle{\sum _{k\in {\bf Z}}}
\exp \left (
\pi i \tau k^2 +2\pi i
xk \right ),
\\ \\
\theta _4(x|\tau)=\displaystyle{\sum _{k\in {\bf Z}}}
\exp \left (
\pi i \tau k^2 +2\pi i
(x+\frac{1}{2})k\right ).
\end{array}
\label{theta}
\eeq
They can be represented as infinite products:
\beq
\label{infp}
\begin{array}{l}
\theta _1(x|\tau)=\displaystyle{
2\sin (\pi x)e^{\pi i\tau /4}
\prod_{k=1}^{\infty}\Bigl ( 1-e^{2\pi i k\tau}\Bigr )
\Bigl ( 1-e^{2\pi i (k\tau +x)}\Bigr )
\Bigl ( 1-e^{2\pi i (k\tau -x)}\Bigr )}\,,
\\ \\
\theta _2(x|\tau)=2\cos (\pi x)e^{\pi i\tau /4}
\displaystyle{
\prod_{k=1}^{\infty}\Bigl ( 1-e^{2\pi i k\tau}\Bigr )
\Bigl ( 1+e^{2\pi i (k\tau +x)}\Bigr )
\Bigl ( 1+e^{2\pi i (k\tau -x)}\Bigr )}\,,
\\ \\
\theta _3(x|\tau)=\displaystyle{
\prod_{k=1}^{\infty}\Bigl ( 1-e^{2\pi i k\tau}\Bigr )
\Bigl ( 1+e^{\pi i (2k\tau -\tau +2x)}\Bigr )
\Bigl ( 1+e^{\pi i (2k\tau -\tau -2x)}\Bigr )}\,,
\\ \\
\theta _4(x|\tau)=\displaystyle{
\prod_{k=1}^{\infty}\Bigl ( 1-e^{2\pi i k\tau}\Bigr )
\Bigl ( 1-e^{\pi i (2k\tau -\tau +2x)}\Bigr )
\Bigl ( 1-e^{\pi i (2k\tau -\tau -2x)}\Bigr )}\,.
\end{array}
\eeq
The $\wp$-function is related to the $\theta_1$ as follows:
\beq
\label{wei}
\wp (x|1/2, \tau /2)=-\frac{d^2}{dx^2}\,\mbox{log}\,
\theta_1(x|\tau )+\mbox{const}\,.
\eeq
Throughout the paper we write
$\theta _a(x|\tau)=\theta _a(x)$.

The
transformation properties of the
theta-functions for shifts by (half) periods are:
\beq
\label{periods}
\theta_a (x\pm 1)=(-1)^{\delta _{a,1}+\delta _{a,2}}
\theta_a (x)\,,
\;\;\;\;\;
\theta_a (x\pm \tau )=(-1)^{\delta _{a,1}+\delta _{a,4}}
e^{-\pi i \tau \mp 2\pi i x}
\theta_a (x)\,,
\eeq
\beq
\label{half}
\begin{array}{l}
\theta_1 (x\pm \frac{1}{2})=\pm
\theta_2 (x)\,,\\ \\
\theta_1 (x\pm \frac{\tau}{2})=\pm i
e^{-\frac{1}{4}\pi i \tau \mp \pi i x}
\theta_4 (x)\,, \\ \\
\theta_1 (x\pm \frac{1+\tau}{2})=\pm
e^{-\frac{1}{4}\pi i \tau \mp \pi i x}
\theta_3 (x)\,.
\end{array}
\eeq

\end{document}